\documentclass{article}
\usepackage{amsthm}
\usepackage{amsmath}
\usepackage{amssymb}

\title{Non-homogeneous locally free actions\\ of the affine group}
\author{Masayuki ASAOKA
\footnote{Partially supported by JSPS PostDoctoral Fellowships for
 Research Abroad
 and Grant-in-Aid for Encouragement of Young
 Scientists (B) (No.19740085)}\\
Department of Mathematics, Kyoto University
}

\def\RR{\mathbb{R}}

\def\HH{\mathbb{H}}

\def\cO{{\cal O}}
\def\cF{{\cal F}}

\def\GA{{\mbox{\rm{GA}}}}
\def\AGamma{{{\cal A}_\Gamma}}

\def\PSL{{\rm PSL}}

\def\wSL{{\widetilde{\PSL}}}
\def\psl{{\mathfrak{sl}_2(\RR)}}
\def\Solv{{\it Solv}}

\def\Per{{\rm Per}}

\def\bsl{{\backslash}}
\def\ra{{\rightarrow}}
\def\st{{\;|\;}}

\def\vphi{\varphi}
\def\cvphi{\check{\varphi}}
\def\del{\partial}

\def\chM{\check{M}}

\def\chF{\check{\cF}}

\def\chPhi{\check{\Phi}}

\def\hsp{{\hspace{3mm}}}

\newcommand\ch[1]{{\check{#1}}}
\newcommand\Pair[2]{{\langle {#1},{#2} \rangle}}

\theoremstyle{plain}
\newtheorem{thm}{Theorem}[section]
\newtheorem{prop}[thm]{Proposition}
\newtheorem{lemma}[thm]{Lemma}
\newtheorem{cor}[thm]{Corollary}

\newtheorem*{Mthm}{Main Theorem}

\newtheorem*{Livthm}{The Livschitz Theorem}

\theoremstyle{definition}

\newtheorem{prob}[thm]{Problem}

\theoremstyle{remark}

\begin{document}
\maketitle

\noindent
{\bf Abstract:}
We classify smooth locally free actions of the real affine group
 on closed orientable three-dimensional manifolds up to smooth conjugacy.
As a corollary, there exists a non-homogeneous action
 when the manifold is the unit tangent bundle of
 a closed surface with a hyperbolic metric.
\medskip



\section{Introduction}
\label{sec:introductioon}
%
%
Let $\GA$ be the group of orientation preserving
 affine transformations of the real line. 
In 1979, Ghys \cite{Gh0} showed that any volume-preserving
 locally free $\GA$-action on a closed three-dimensional manifold
 is smoothly conjugate to a homogeneous action.
It is natural to ask whether  the assumption
 on an invariant volume is removable or not,
 and it was open for almost thirty years.
In this paper, we give a classification
 of locally free $\GA$-actions on closed orientable three-dimensional
 manifolds up to smooth conjugacy.
As a corollary,
 there {\it exists} an action with no invariant volume.

Let $G$ be a Lie group and $M$ a manifold.
We say that a $C^\infty$ right $G$-action $\rho:M \times G \ra M$
 is {\it locally free} if the isotropy subgroup
 $\{g \in G \st \rho(p,g)=p\}$ is discrete in $G$
 for any $p \in M$.
By $\cO_\rho(p)$, we denote the orbit $\{\rho(p,g) \st g \in G\}$
 of $p \in M$.
If the action $\rho$ is locally free,
 the decomposition $\cO_\rho=\{\cO_\rho(p) \st p \in M\}$ 
 forms a $C^\infty$ foliation on $M$.
We say that two right $G$-actions $\rho_1$ and $\rho_2$
 on manifolds $M_1$ and $M_2$
 are {\it $C^\infty$-conjugate} if there exists an isomorphism
 $\sigma$ of $G$ and a $C^\infty$ diffeomorphism $H$ from $M_1$ to $M_2$
 such that $H(\rho_1(p,g))=\rho_2(H(p), \sigma(g))$
 for any $g \in G$ and $p \in M_1$.
The map $H$ is called {\it a $C^\infty$ conjugacy}
 between $\rho_1$ and $\rho_2$.

There are two classical examples of locally free $\GA$-actions
 on closed three-dimensional manifolds.
Let $\PSL(2,\RR)$ be the group of orientation preserving
 projective transformations of the real projective line $\RR P^1$
 and $\wSL(2,\RR)$ its universal covering group.
The group $\wSL(2,\RR)$
 contains a closed subgroup $H$ which is isomorphic to $\GA$,
 {\it e.g.}, the group of elements which fix
 $\infty \in \RR P^1 =\RR \cup \{\infty\}$.
Fix an isomorphism $i:\GA \ra H$.
For each cocompact lattice $\Gamma$,
 {\it the standard $\GA$-action} $\rho_\Gamma$
 on $\Gamma \backslash \wSL(2,\RR)$ is given
 by $\rho_\Gamma(\Gamma g,h)=\Gamma(g \cdot i(h))$.
A similar construction gives
 another locally free action $\rho_{solv}$ when we start with
 a three-dimensional solvable subgroup
\begin{equation*}
\Solv=\{h_{t,x,y}:(u,v) \mapsto (e^{-t}u+x,e^t v+y) \st t,x,y \in \RR\} 
\end{equation*}
 of the group of affine transformations of $\RR^2$.
In this case, $\GA$ is identified
 with the subgroup $\{h_{t,x,0} \st t,x \in \RR\}$.
We say that a $\GA$-action is {\it homogeneous}
 if it is $C^\infty$-conjugate to one of the above examples.

As mentioned at the beginning,
 Ghys showed the rigidity of volume-preserving actions
 in his doctoral thesis \cite{Gh0}.
\begin{thm}
[\cite{Gh0}, see also \cite{Gh}]
\label{thm:Ghys 1}
Let $M$ be a closed, connected, oriented, and three-dimensional manifold.
If a $C^\infty$ locally free action of $\GA$ on $M$
 admits a continuous invariant volume form, then it is homogeneous.
\end{thm}
The following problem is quite natural, but was open for thirty years.
\begin{prob}
[{\cite[p.525]{Gh}},
 see also {\cite[p.2]{Be2}, \cite[p.1841]{MM}}]
\label{prob:rigidity}
Remove the assumption on an invariant volume,
 or construct a $C^\infty$ locally free action of $\GA$
 on a closed three-dimensional manifold which is not homogeneous.
\end{prob}
There are several partial solutions to the problem.
We say that two actions $\rho_1$ and $\rho_2$ on manifolds $M_1$ and $M_2$
 are {\it $C^\infty$ orbit equivalent}
 if there exists a $C^\infty$ diffeomorphism $H$ from $M_1$ to $M_2$
 such that $H(\cO_{\rho_1}(p))=\cO_{\rho_2}(H(p))$
 for any $p \in M_1$.

\begin{thm}
[\cite{Gh}]
\label{thm:Ghys 3}
Any $C^\infty$ locally free $\GA$-action
 on a three-dimensional rational homology sphere
 preserves a continuous volume.
In particular, it is a homogeneous action.
\end{thm}

\begin{thm}
[\cite{Gh2}]
\label{thm:Ghys 4}
Let $\Gamma$ be a cocompact lattice of $\Solv$.
Then, any $C^\infty$ locally free action of $\GA$ on 
 $\Gamma \backslash \Solv$ is homogeneous.
\end{thm}

\begin{thm}
[\cite{As1}]
\label{thm:Asaoka 1}
Let $M$ be a closed, connected, oriented, and three-dimensional manifold
 and $\rho$ a $C^\infty$ locally free action of $\GA$ on $M$.
Then, $\rho$ is $C^\infty$ orbit equivalent to a homogeneous action.
In particular,
 $M$ is diffeomorphic to
 $\Gamma \backslash \wSL(2,\RR)$ or $\Gamma \backslash \Solv$
 with a cocompact lattice $\Gamma$.
\end{thm}
We will give a proof of Theorem \ref{thm:Ghys 4}
 in Appendix since it has not been published.

Let $\Gamma$ be a cocompact lattice of $\wSL(2,\RR)$.
We denote the quotient space $\Gamma \bsl \wSL(2,\RR)$ by $M_\Gamma$.
Let $\AGamma$ be the set of $C^\infty$ locally free actions
 of $\GA$ on $\Gamma \bsl \wSL(2,\RR)$
 whose orbit foliation coincides with that of the standard $\GA$-action.
By the above results,
 it is sufficient to classify
 actions in $\AGamma$ up to $C^\infty$-conjugacy
 for the complete classification of
 locally free  $\GA$-actions on closed, orientable,
 and three-dimensional manifold.\footnote{
We can also classify actions on non-orientable manifolds
 by taking an orientable double cover.}
\begin{Mthm}
There exists an open subset $\Delta_\Gamma$
 of $H^1(M_\Gamma,\RR)$
 and a family $(\rho_a)_{a \in \Delta_\Gamma}$ of actions $\AGamma$
 such that
\begin{enumerate}
\item any $\rho \in \AGamma$
 is $C^\infty$-conjugate to $\rho_a$ for some $a \in \Delta_\Gamma$,
 and
\item $a=a'$ if and only if
 $\rho_a$ is $C^\infty$-conjugate to $\rho_{a'}$
 by a diffeomorphism homotopic to the identity.
\end{enumerate}
In particular,
 if $M_\Gamma$ is not a rational homology sphere,
 then $\AGamma$ contains non-homogeneous actions.
\end{Mthm}

We will not prove the smoothness
 of the family $(\rho_a)_{\alpha \in \Delta_\Gamma}$ in this paper.
It need more efforts, and will be shown in a forthcoming paper.
\medskip

The organization of the paper is the following:
In Section \ref{sec:non-homogeneous},
 we associate a $\GA$-action $\ch{\rho}_\omega$ 
 with any given closed one-form $\omega$ on $\Gamma \bsl \wSL(2,\RR)$
 satisfying a mild condition.
The action is not homogeneous
 if the cohomology class of $\omega$ is non-zero. 
As a corollary, we obtain the existence of non-homogeneous actions
 when $\Gamma \bsl \wSL(2,\RR)$ is not a rational homology sphere.
In Section \ref{sec:classification},
 we define a natural parametrization $\bar{a}_\Gamma$ of actions
 and determine its image $\Delta_\Gamma$.
In fact, we show that
 any action in $\AGamma$
 is $C^\infty$-conjugate to
 an action in the family $(\ch{\rho}_\omega)$
 constructed in Section \ref{sec:non-homogeneous}.
In Appendix,
 we give a proof of Theorem \ref{thm:Ghys 4}
 and discuss the regularity of the unstable foliation
 of the Anosov flow associated with a non-homogeneous action.

\paragraph{Acknowledgements}
The author wrote the first version of this article when he stayed
 at Unit\'e de Math\'ematiques Pures et Appliqu\'ees,
 \'Ecole Normale Sup\'erieure de Lyon.
He would like to thank the members of UMPA,
 especially Professor \'Etienne Ghys,
 for their warm hospitality.
The author is also grateful to Professor Shigenori Matsumoto
 and an anonymous referee for many suggestions to improve the paper.
At last, I would like to thank my wife, Junko for encouragements.

\section{Non-homogeneous actions}
\label{sec:non-homogeneous}
%
%

In this section, we construct a non-homogeneous $\GA$-action.
More precisely, for any given closed one-form $\omega$ on
 $M_\Gamma=\Gamma \bsl \wSL(2,\RR)$ with a mild condition,
 we construct a $C^\omega$ action $\ch{\rho}_\omega$ on
 a manifold diffeomorphic to $M_\Gamma$
 such that $\det D\ch{\rho}_\omega$ is controlled by $\omega$
 (Theorem \ref{thm:action}).
The action $\ch{\rho}_\omega$ is not homogeneous
 if the cohomology class of $\omega$ is non-zero
 (Proposition \ref{prop:action}).
As a corollary,
 if $H^1(M_\Gamma)$ is non-trivial,
 then $M_\Gamma$ admits a $C^\omega$ non-homogeneous action.

The main ingredient of the construction is
 deformation of a codimension-one Anosov system
 along the strong stable foliation.
It is invented by E.Cawley \cite{Ca}
 in order to describe the moduli space of
 two dimensional Anosov diffeomorphisms.
Here, we apply her technique to codimension-one Anosov flows
 with a $C^\omega$ invariant splitting.
It produces another Anosov flow with
 the $C^\omega$ strong stable foliation
 whose rate of contraction is constant.
Such a flow induces a locally free $\GA$-action naturally.

In Subsection \ref{sec:Margulis}, we recall the definitions
 of an Anosov flow and the Margulis measure associated with it.
Deformation of Anosov flow along the strong stable foliation
 is done in Subsection \ref{sec:deform}.
Finally, we construct non-homogeneous actions
 associated with a closed one-form in Subsection \ref{sec:action}

\subsection{Anosov flows and the Margulis measure}
\label{sec:Margulis}
Let $\Phi$ be a $C^1$ flow on a closed manifold $M$
 without stationary points.
Its differential defines a flow $D\Phi$
 on the tangent bundle $TM$ of $M$.
Let $T\Phi$ be the one-dimensional subbundle of $TM$
 which is tangent to the orbits of $\Phi$.
The flow $\Phi$ is called {\it Anosov} 
 if there exists a continuous $D\Phi$-invariant splitting
 $TM=T\Phi \oplus E^{ss} \oplus E^{uu}$
 and $T>0$ such that
\begin{equation*}
 \max\left\{
  \|D\Phi^T|_{E^{ss}(p)}\|,
  \|D\Phi^{-T}|_{E^{uu}(p)}\| \right\} < 1/2
\end{equation*}
 for any $p \in M$,
 where $\|\cdot \|$ is the operator norm with respect to
 a Riemannian metric of $M$.
The above splitting is unique for $\Phi$
 and is called {\it the Anosov splitting} of $\Phi$.
The $D\Phi$-invariant subbundles
 $E^{ss}$, $E^{uu}$, $E^s=T\Phi \oplus E^{ss}$,
 and $E^u=T\Phi \oplus E^{uu}$ of $TM$
 generate continuous foliations
 $\cF^{ss}$, $\cF^{uu}$, $\cF^s$, and $\cF^u$.
They are called
 {\it the strong stable, strong unstable,
 weak stable}, and {\it weak unstable foliations}, respectively.

By $\Per(\Phi)$, we denote the union of periodic orbits of $\Phi$.
For $p \in \Per(\Phi)$,
 let $\tau(p;\Phi)$ be the (minimal) period of $p$.
We define a loop $\gamma(p;\Phi):[0,1] \ra M$
 by $\gamma(p;\Phi)(t)=\Phi^{t \cdot \tau(p;\Phi)}(p)$.
Put
\begin{align*}
J(p;\Phi)
 & =\log \left|\det D\Phi^{\tau(p;\Phi)}_p\right|,\\
J^s(p;\Phi)
 & =\log \left|\det \left(D\Phi^{\tau(p;\Phi)}|_{E^{ss}(p)}\right)\right|,\\
J^u(p;\Phi)
 & =\log \left|\det \left(D\Phi^{\tau(p;\Phi)}|_{E^{uu}(p)}\right)\right|.
\end{align*}
Remark that $J$, $J^s$ and $J^u$ does not depend on the
 choice of the Riemannian metric.

When two points $p,q \in M$ are sufficiently close to each other,
 we can define the holonomy map $\tilde{h}_{pq}$ of $\cF^{uu}$
 from a neighborhood $\tilde{U}_{pq}$ of $p$ in $\cF^s(p)$ to $\cF^s(q)$.
We say that a flow is {\it topologically transitive} if it admits a dense orbit.
\begin{thm}
[Margulis \cite{Ma}]
If $\Phi$ is a topologically transitive Anosov flow,
 then there exists a family $(\tilde{\mu}_p)_{p \in M}$
 which has the following properties:
\begin{enumerate}
\item $\tilde{\mu}_p$ is a locally finite Borel measure on $\cF^s(p)$
 which is non-atomic and positive on each non-empty open subset of $\cF^s(p)$.
\item  $\tilde{\mu}_p=\tilde{\mu}_q$ if $\cF^s(p)=\cF^s(q)$.
\item The family $(\tilde{\mu}_p)_{p \in M}$ is holonomy invariant,
 {\it i.e.}, if $p$ and $q$ are sufficiently close to each other,
 then $\tilde{\mu}_q(\tilde{h}_{pq}(A))=\tilde{\mu}_p(A)$
 for any Borel subset $A$ of $\tilde{U}_{pq}$.
\item There exists $\lambda>0$ such that
\begin{equation*}
 \tilde{\mu}_p(\Phi^t(A))=e^{-\lambda t}\cdot \tilde{\mu}_p(A)
\end{equation*}
 for any $p \in M$, $t \in \RR$, and any Borel subset $A$ of $\cF^s(p)$.
\end{enumerate}
\end{thm}

As mentioned in Chapter 3 of \cite{Ma}, 
 the family $(\tilde{\mu}_p)_{p \in M}$ induces
 a family $(\mu_p)_{p \in M}$ of measures on leaves of $\cF^{ss}$
 such that
\begin{equation*}
 \tilde{\mu}_p(\{\Phi^t(q) \st (q,t) \in \tilde{A} \times (0,\epsilon)\})
 =  (1-e^{-\lambda \epsilon}) \cdot \mu_p(\tilde{A}).
\end{equation*}
 for any small Borel set $\tilde{A} \subset \cF^{ss}(p)$ 
 and any small $\epsilon>0$.
We call the family $(\mu_p)_{p \in M}$
 {\it the Margulis measure} of $\Phi$ along $\cF^{ss}$.
It satisfies the following properties:
\begin{itemize}
\item
 For any small Borel subset $\tilde{A} \subset \cF^s(p)$
 of the form
\begin{equation*}
 \tilde{A}=\{\Phi^t(p') \st (p',t) \in A \times [\eta_-(p'), \eta_+(p')]\},
\end{equation*}
 where $A$ is a Borel subset of $\cF^{ss}(p)$, and
 $\eta_-$ and $\eta_+$ are measurable functions on $A$,
\begin{equation}
\label{eqn:mu-p 1}
 \tilde{\mu}_p(\tilde{A})=\int_A
 \{e^{-\lambda \cdot \eta_-(p')}-e^{-\lambda \cdot \eta_+(p')}\} d\mu_p(p').
\end{equation}
\item For any $p \in M$, $t \in \RR$,
 and any small Borel subset $A$ of $\cF^{ss}(p)$,
\begin{equation}
\label{eqn:mu-p 2}
 \mu_{\Phi^t(p)}(\Phi^t(A))=e^{-\lambda t}\cdot \mu_p(A)
\end{equation}
\end{itemize}

When two points $p$ and $q$ in $M$ are sufficiently close
 to each other,
 we can define the holonomy map $h_{pq}$ of $\cF^u$
 from a neighborhood $U_{pq}$ of $p$ in $\cF^{ss}(p)$
 to $\cF^{ss}(q)$.
\begin{prop}
\label{prop:smooth derivative}
Suppose that the flow $\Phi$
 and the foliations $\cF^{ss}$, $\cF^{uu}$,  $\cF^s$, $\cF^u$
 are real analytic.
Then, for any given $p \in M$, the Radon-Nikodym derivative
 $\frac{d (\mu_q \circ h_{pq})}{d \mu_p}(p')$
 is a $C^\omega$ function of $q \in M$ and $p' \in U_{pq}$.
\end{prop}
\begin{proof}
Since $\cF^s$ is subfoliated by $\cF^{ss}$,
 there exists a neighborhood $U$ of $p$ in $M$
 and a $C^\omega$ function $\eta$ on $U \times U$ such that
 $\eta(p',p')=0$ and
 $\Phi^{\eta(p',q)}(\tilde{h}_{p'q}(p'))=h_{p'q}(p')$
 for any $(p',q) \in U \times U$.

Take a small Borel set $A \subset U \cap U_{pq}$
 and a small positive number $\epsilon>0$.
Since $\Phi$ preserves the foliation $\cF^{uu}$,
 we have $\tilde{h}_{pq} \circ \Phi^t=\Phi^t \circ \tilde{h}_{pq}$
 for any $t \in [0,\epsilon]$.
It implies
\begin{align*}
(1-e^{-\lambda \epsilon}) \mu_q(h_{pq}(A))
 & =\tilde{\mu}_q(\{\Phi^{t+\eta(p',q)} \circ \tilde{h}_{pq}(p')
  \st (p',t) \in A \times (0,\epsilon)\})\\
 & = (\tilde{\mu}_q \circ \tilde{h}_{pq})
 (\{\Phi^{t+\eta(p',q)}(p') \st (p',t) \in A \times (0,\epsilon)\})\\
 & = \tilde{\mu}_p
 (\{\Phi^{t+\eta(p',q)}(p') \st (p',t) \in A \times (0,\epsilon)\})\\
 & = \int_A \{e^{-\lambda \cdot \eta(p',q)}-e^{-\lambda(\eta(p',q)+\epsilon)}\}
 d \mu_p(p').
\end{align*}
In particular, 
\begin{equation*}
(\mu_q \circ h_{pq})(A)=\int_A e^{-\lambda \eta(p',q)} d \mu_p(p'),
\end{equation*}
 and hence,
\begin{equation*}
\frac{d (\mu_q \circ h_{pq})}{d \mu_p} (p')
 =\exp(-\lambda \cdot \eta(p',q))
\end{equation*}
 for any $p' \in \cF^{ss}(p)$ sufficiently close to $p$.
The right-hand side is a $C^\omega$ function of $(p',q)$.
\end{proof}

\subsection{Deformation of codimension-one Anosov flows}
\label{sec:deform}
In this subsection, we deform an Anosov flow 
 into another Anosov flow with constant contraction.
\begin{prop}
\label{prop:construction} 
Let $M$ be an $n$-dimensional closed manifold
 and $\Phi$ a $C^\omega$ topologically transitive Anosov flow on $M$.
Suppose that  $\Phi$ admits the $C^\omega$ Anosov splitting
 and the strong stable foliation is one-dimensional and orientable.
Then, there exists
 another $C^\omega$ Anosov flow $\chPhi$ on a $C^\omega$ closed manifold $\chM$,
 a $C^\omega$ vector field $Y^s$ on $\chM$,
 a homeomorphism $H:M \ra \chM$,
 and a constant $\lambda>0$ which satisfy the following properties:
\begin{itemize}
\item $H \circ \Phi^t= \chPhi^t \circ H$ for any $t \in \RR$.
\item $J^u(H(p);\ch{\Phi})=J^u(p;\Phi)$ for any $p \in \Per(\Phi)$.
\item  $Y^s$ generates the strong stable foliation of $\chPhi$
 and $D\chPhi^t(Y^s) =e^{-\lambda t} \cdot Y^s$ for any $t \in \RR$.
In particular,
 $J^s(\ch{p};\ch{\Phi}) = -\lambda \cdot \tau(\ch{p};\ch{\Phi})$
 for any $p \in \Per(\ch{\Phi})$.
\end{itemize}
\end{prop}
The rest of this subsection is devoted to the proof.
We follow Cawley's idea in \cite{Ca},
 where she deformed two-dimensional Anosov diffeomorphisms.
Although the deformed diffeomorphism is of class $C^1$ in her case,
 we obtain a $C^\omega$ flow in our case
 because of the regularity of the invariant splitting
 and the constant contraction along $\cF^{ss}$
 with respect to the Margulis measure.

Let $\cF^{ss}$ and $\cF^u$ be the strong stable foliation
 and the weak unstable foliation of $\Phi$.
We denote the coordinate system of $\RR^n$
 by $(x_1,\cdots,x_{n-1},y)$.
For each $p \in M$, we take an open neighborhood $U_p$ of $p$
 and a $C^\omega$ diffeomorphism $\vphi_p:U_p \ra (-1,1)^n$
 such that $\vphi_p(p)=(0,0)$ and
\begin{gather*}
\vphi_p^{-1}(\{x\} \times (-1,1)) \subset \cF^{ss}(\vphi_p^{-1}(x,y)),\\
\vphi_p^{-1}((-1,1)^{n-1} \times \{y\}) \subset \cF^u(\vphi_p^{-1}(x,y))
\end{gather*}
 for any $(x,y) \in \RR^{n-1} \times \RR$.
Let $\pi_p:U_p \ra \vphi_p^{-1}((-1,1)^{n-1} \times 0)$
 be the projection along $\cF^{ss}$.
For $q \in U_p$, let $h_{pq}:(\cF^{ss}|_{U_p})(p) \ra (\cF^{ss}|_{U_p})(q)$
 be the holonomy map along $\cF^u$.
They satisfy
\begin{gather*}
\pi_p \circ \vphi_p^{-1}(x,y)=\vphi_p^{-1}(x,0),\\
h_{pq} \circ \vphi_p^{-1}(0,y)=\vphi_p^{-1}(x_q,y) 
\end{gather*}
 for any $(x,y) \in (-1,1)^{n-1} \times (-1,1)$,
 where $\vphi_p(q)=(x_q, y_q)$.

Fix an orientation of $\cF^{ss}$.
Let $(\mu_p)_{p \in M}$ be the Margulis measure of $\Phi$
 along $\cF^{ss}$.
For an oriented interval $I \subset \cF^{ss}(p)$,
we put $\nu(I)=\mu_p(I)$ if
 the orientations of $I$ and $\cF^{ss}(p)$ coincide
 or $\nu(I)=-\mu_p(I)$ otherwise.
For $p \in M$ and $q_-,q_+ \in \cF^{ss}(p)$,
 let $[q_-,q_+]$ be the oriented interval in $\cF^{ss}(p)$
 which connects $q_-$ and $q_+$
 and whose positive end point is $q_+$.

\begin{lemma}
\label{lemma:analytic holonomy} 
For any given $y_1,y_2 \in (-1,1)$, the function
\begin{equation*}
x \mapsto \nu([\vphi_p^{-1}(x,y_1),\vphi_p^{-1}(x,y_2)]) 
\end{equation*}
 from $(-1,1)^{n-1}$ to $\RR$ is real-analytic.
\end{lemma}
\begin{proof}
Without loss of generality, we may assume that $y_1<y_2$.
Put $q(x)=\vphi_p^{-1}(x,0)$ for $x \in (-1,1)^{n-1}$.
Then,
\begin{align*}
|\nu([\vphi_p^{-1}(x,y_1),\vphi_p^{-1}(x,y_2)])|
 & =\mu_{q(x)}(\vphi_p^{-1}(x \times [y_1,y_2]))\\
 & =(\mu_{q(x)} \circ h_{pq(x)})(\vphi_p^{-1}(0 \times [y_1,y_2]))\\
 & = \int_{\vphi_p^{-1}(0 \times [y_1,y_2])}
  \frac{d(\mu_{q(x)} \circ h_{pq(x)})}{d\mu_p} \, d\mu_p.
\end{align*}
By Proposition \ref{prop:smooth derivative},
 the last term is real-analytic with respect to $x$.
\end{proof}

We define a function $\ch{y}_p:U_p \ra \RR$
 and a map $\cvphi_p:U_p \ra \RR^n$ by
\begin{gather*}
 \ch{y}_p(q)=\nu([\pi_p(q),q]),\\
 \cvphi_p(q)=(x,\ch{y}_p(q))
\end{gather*}
 for $q \in U_p$ with $\vphi_p(q)=(x,y)$.
Since $\mu_p$ is locally finite, non-atomic,
 and positive on any non-empty open subset of $\cF^{ss}(p)$,
 the map $\cvphi_p$ is a homeomorphism onto
 an open subset of $\RR^n$.
Remark that
\begin{equation*}
\pi_q \circ \cvphi_p^{-1}(x,y) = \vphi_p^{-1}(x,y_q)
\end{equation*}
 for any $q \in U_p$ with $\vphi_p(q)=(x_q,y_q)$
 and any $(x,y) \in \cvphi_p(U_p \cap U_q)$.
\begin{lemma}
\label{lemma:smooth coordinate} 
For $p \in M$ and $q \in U_p$,
 the map $\cvphi_q \circ \cvphi_p^{-1}$ is real-analytic.
The differential $D(\cvphi_q \circ \cvphi_p^{-1})$
 preserves the vector field $(\del/\del y)$.
\end{lemma}
\begin{proof}
Since $\vphi_p$ and $\vphi_q$ are adapted
 to the pair $(\cF^u, \cF^{ss})$ of foliations,
 there exist $C^\omega$ coordinate changes $\psi$ and $\psi'$
 such that $\vphi_q \circ \vphi_p^{-1}(x,y)=(\psi(x),\psi'(y))$.
For $p' \in U_p \cap U_q$ with $(x,y)=\cvphi_p(p')$,
 we have
\begin{align*}
\cvphi_q \circ \cvphi_p^{-1}(x,y)
 & = (\psi(x),\nu([\pi_q(p'),p']))\\
 & = (\psi(x), \nu([\pi_p(p'),p'])-\nu([\pi_p(p'),\pi_q(p')]))\\
 & = (\psi(x), y-\nu([\vphi_p^{-1}(x,0),\vphi_p^{-1}(x,y_q)])),
\end{align*}
 where $\vphi_p(q)=(x_q,y_q)$.
By Lemma \ref{lemma:analytic holonomy},
 the last term is real-analytic with respect to $(x,y)$.
The above equation also implies that
 $D(\cvphi_q \circ \cvphi_p^{-1})$ preserves $(\del/\del y)$.
\end{proof}

The lemma implies that $(\cvphi_p)_{p \in M}$ defines
 a $C^\omega$-structure on $M$.
By $\ch{M}$, we denote the manifold $M$ endowed
 with this $C^\omega$-structure.
Let $H:M \ra \ch{M}$ be the identity map as a set.
It is a homeomorphism but is not a diffeomorphism in general
 since $\cvphi_p \circ \vphi_p^{-1}$ is continuous
 but may not be smooth.
We define a continuous foliation $\chF^u$ on $\chM$
 by $\chF^u(H(p))=H(\cF^u(p))$.
\begin{lemma}
Each leaf of $\chF^u$ is a $C^\omega$ immersed manifold
 and the restriction of $H$ to a leaf
 is a $C^\omega$ diffeomorphism.
\end{lemma}
\begin{proof}
For $p \in M$,
 $\vphi_p((-1,1)^{n-1} \times 0)$ is a neighborhood of $p$
 in $\cF^u(p)$
 and $\cvphi_p((-1,1)^{n-1} \times 0)$
 is a neighborhood of $H(p)$ in $\chF^u(H(p))$.
The latter implies that $\chF^u(H(p))$ is a $C^\omega$ immersed manifold.
Since $\cvphi_p \circ \vphi_p^{-1}(x,0)=(x,0)$,
 the restriction of $H$ to $\cF^u(p)$ is a $C^\omega$ diffeomorphism.
\end{proof}

Let $TM=T\Phi \oplus E^{ss} \oplus E^{uu}$
 be the Anosov splitting of $\Phi$.
By the above lemma,
 the subbundle $\ch{E}^{uu}=DH(E^{uu})$ is well-defined.
Since $D(\cvphi_q \circ \cvphi_p^{-1})$ preserves
 the vector field $\del /\del y$,
 we can define a $C^\omega$ vector field $Y^s$ on $\chM$ by
 $Y^s=(D\cvphi_p)^{-1}(\del/\del y)$ on $U_p$.
We also define a flow $\ch{\Phi}$ on $\chM$ by
 $\ch{\Phi}^t = H \circ \Phi^t \circ H^{-1}$.
Recall that there exists $\lambda>0$ such that
 $\mu_{\Phi^t(p)} \circ \Phi^t=e^{-\lambda t} \mu_p$
 for any $p \in M$ and $t \in \RR$.

Now, we check that the quadruple $(\ch{\Phi},H,Y^s,\lambda)$ satisfies
 the conditions in Proposition \ref{prop:construction}.
\begin{proof}
[Proof of Proposition \ref{prop:construction}]
Fix $p \in M$.
Since $\Phi$ preserves the foliations $\cF^u$ and $\cF^{ss}$,
 there exists a $C^\omega$ local flow $\Psi_p$ on $(-1,1)^{n-1}$,
 a neighborhood $V_p \subset U_p$ of $p$,
 and $\epsilon>0$ such that
 $\vphi_p \circ \Phi^t \circ \vphi_p^{-1}(x,y)=(\Psi_p^t(x),y)$
 for any $(x,y) \in \vphi_p(V_p)$ and $t \in (-\epsilon,\epsilon)$.
Take $q \in V_p$ and $t \in (-\epsilon,\epsilon)$.
Put $(x,y)=\cvphi_p(q)$ and $(x',y')=\cvphi_p \circ \ch{\Phi}^t(q)$.
Then, $x'=\Psi_p^t(x)$ and
\begin{align*}
y' & = \nu([\pi_p(\Phi^t(q)),\Phi^t(q)])
 =\nu(\Phi^t([\pi_p(q),q]))
 =e^{-\lambda t} \nu([\pi_p(q),q])
  = e^{-\lambda t} y.
\end{align*}
Therefore,
\begin{equation*}
 \cvphi_p \circ \ch{\Phi}^t \circ \cvphi_p^{-1}(x,y)
 =(\Psi_p^t(x), e^{-\lambda t} y).
\end{equation*}
It implies that $\ch{\Phi}$ is a $C^\omega$ flow
 and $D\ch{\Phi}^t(Y^s)=e^{-\lambda t} \cdot Y^s$.
Since $DH|_{E^u}$ is well-defined
 and $DH(E^{uu})=\ch{E}^{uu}$,
 the flow $D\ch{\Phi}$ preserves $\ch{E}^{uu}$ and 
 there exists a $C^\infty$ norm $\|\cdot\|$ on $T\ch{M}$
 and $T>0$ such that
 $\|D\ch{\Phi}^{-T}|_{\ch{E}^{uu}(p)}\|<1/2$ for any $p \in \ch{M}$.
Hence, $\ch{\Phi}$ is an Anosov flow
 and its Anosov splitting
 is $T\chM=T\ch{\Phi} \oplus \RR Y^s \oplus \ch{E}^{uu}$.
Since $DH|_{E^{uu}}$ is well-defined,
 $J^u(H(q);\ch{\Phi})=J^u(q;\Phi)$ for any $q \in \Per(\Phi)$.
\end{proof}

\subsection{Construction of a non-homogeneous action}
\label{sec:action}
Fix a cocompact lattice $\Gamma$ of $\wSL(2,\RR)$.
We denote the quotient space $\Gamma \bsl \wSL(2,\RR)$ by $M_\Gamma$
 and the Lie algebra of $\wSL(2,\RR)$ by $\psl$.
The Lie algebra $\psl$ is naturally
 identified with the Lie algebra of trace-free real $2\times 2$-matrices.
It is generated by three elements
\begin{equation*}
X=
\begin{bmatrix}
1/2 & 0 \\ 0 & -1/2
\end{bmatrix},\hsp
S=
\begin{bmatrix}
0 & 1 \\ 0 & 0
\end{bmatrix},\hsp
U=
\begin{bmatrix}
0 & 0 \\ 1 & 0
\end{bmatrix}.
\end{equation*}
Let $(X^t)_{t \in \RR}$, $(S^t)_{t \in \RR}$, and $(U^t)_{t \in \RR}$
 be the one-parameter subgroups of $\wSL(2,\RR)$ corresponding to
 $X$, $S$, and $U$, respectively.
Since $[X,S]=S$,
 the Lie subalgebra spanned by $X$ and $S$
 generates a subgroup of $\wSL(2,\RR)$ which is isomorphic to $\GA$.
{\it In the rest of the paper, we identify this group with $\GA$.}
Under the identification,
 the standard $\GA$-action $\rho_\Gamma$ on $M_\Gamma$
 is written as $\rho_\Gamma(\Gamma g,h)=\Gamma(gh)$.

The elements $X$, $S$, and $U$ of $\psl$ can be
 identified with the left invariant vector fields on $\wSL(2,\RR)$.
They naturally induce
 vector fields $X_\Gamma$, $S_\Gamma$, and $U_\Gamma$
 on $M_\Gamma$.
We define a flow $\Phi_\Gamma$ on $M_\Gamma$ by
 $\Phi_\Gamma^t(\Gamma g)=\Gamma(g X^t)$.
It is generated by $X_\Gamma$
 and satisfies $D\Phi_\Gamma^t(S_\Gamma) = e^{-t} \cdot S_\Gamma$
 and $D\Phi_\Gamma^t(U_\Gamma) = e^{t} \cdot U_\Gamma$
 for any $t \in \RR$.
Hence, $\Phi_\Gamma$ is an Anosov flow with the Anosov splitting
 $TM_\Gamma=T\Phi \oplus \RR S_\Gamma \oplus \RR U_\Gamma$ and 
\begin{equation*}
\tau(p;\Phi_\Gamma)=J^u(p;\Phi_\Gamma)=-J^s(p;\Phi_\Gamma) 
\end{equation*}
 for any $p \in \Per(\Phi_\Gamma)$.
For a loop $\gamma$ in $M_\Gamma$, we denote its homology class
 by $[\gamma]$.
The following properties of $\Phi_\Gamma$ are well-known and easy to prove:
\begin{itemize}
\item $\Phi_\Gamma$ is topologically transitive.
\item There exists $p_0 \in \Per(\Phi_\Gamma)$
 such that $[\gamma(p_0;\Phi_\Gamma)]=0$.
\item The set $\{[\gamma(p;\Phi_\Gamma)] \st p \in \Per(\Phi_\Gamma)\}$
 spans $H_1(M_\Gamma,\RR)$.
\item For any $p \in \Per(\Phi_\Gamma)$,
 there exists $p' \in \Per(\Phi_\Gamma)$
 such that $\tau(p,\Phi_\Gamma)=\tau(p',\Phi_\Gamma)$
 and $[\gamma(p';\Phi_\Gamma)]=-[\gamma(p;\Phi_\Gamma)]$.
\end{itemize}

We say that two flows $\Psi_1$ and $\Psi_2$ on manifolds $M_1$
 and $M_2$ are {\it topologically equivalent}
 if there exists a homeomorphism $H:M_1 \ra M_2$
 which sends orbits of $\Psi_1$ to those of $\Psi_2$
 and preserves the orientation of orbits.
The map $H$ is called {\it a topological equivalence} between
 $\Psi_1$ and $\Psi_2$.

Now, we are ready to construct a $\GA$-action associated with
 a closed one-form.
\begin{thm}
\label{thm:action}
Let $\omega$ be a $C^\omega$ closed one-form on $M_\Gamma$
 with $1+\omega(X_\Gamma)>0$.
Then,
 there exists a $C^\omega$ manifold $M_\omega$,
 a $C^\omega$ locally free $\GA$-action $\ch{\rho}_\omega$ on $M_\omega$,
 a homeomorphism $H:M_\Gamma \ra M_\omega$, and
 a constant $\lambda>0$ such that
\begin{enumerate}
\item the flow $\Phi_{\ch{\rho}_\omega}$
 defined by $\Phi_{\ch{\rho}_\omega}^t(\ch{p})=\ch{\rho}_\omega(\ch{p},X^t)$
 is an Anosov flow on $M_\omega$,
\item $H$ is a topological equivalence between
 $\Phi_\Gamma$ and $\Phi_{\ch{\rho}_\omega}$, and
 \item  for any $p \in \Per(\Phi_\Gamma)$,
\begin{align*}
 J^u(H(p);\Phi_{\ch{\rho}_\omega})
 &  = \tau(p;\Phi_\Gamma),\\
 - J^s(H(p);\Phi_{\ch{\rho}_\omega})
 & = \tau(H(p);\Phi_{\ch{\rho}_\omega})
  =\lambda \cdot \left\{\tau(p;\Phi_\Gamma)
 +\Pair{\omega}{\gamma(p;\Phi_\Gamma)}\right\}.
\end{align*}
\end{enumerate}
\end{thm}
\begin{proof}
Fix a $C^\omega$ closed one-form $\omega$
 such that $1+\omega(X_\Gamma)>0$.
The vector field $X_\omega=(1+\omega(X_\Gamma))^{-1}X_\Gamma$
 generates an Anosov flow $\Phi_\omega$
 whose orbit is same as $\Phi_\Gamma$.\footnote{The flow $\Phi_\Gamma$
 is essentially same as the flow investigated
 in \cite{Gh5} and \cite{Sa}.}
Then,
\begin{equation*}
\tau(p;\Phi_\omega) =\tau(p;\Phi_\Gamma)+\Pair{\omega}{\gamma(p;\Phi_\Gamma)}
\end{equation*}
 for any $p \in \Per(\Phi_\Gamma)$.
Put $S_\omega=S_\Gamma-\omega(S_\Gamma) \cdot X_\omega$ 
 and $U_\omega=U_\Gamma-\omega(U_\Gamma) \cdot X_\omega$.
By a direct calculation, we have
\begin{align*}
[S_\omega,X_\omega]
 & =(1+\omega(X_\Gamma))^{-1} \cdot
  \left(
  -S_\Gamma
  +\{X_\Gamma\cdot \omega(S_\Gamma)-S_\Gamma \cdot \omega(X_\Gamma)\}X_\omega
  \right).
\end{align*}
Since $d \omega =0$ \footnote{In the construction of $\ch{\rho}_\omega$,
 this is the only place we use $d \omega=0$.},
 it implies that $[S_\omega,X_\omega]=-(1+\omega(X_\Gamma))^{-1}S_\omega$.
Similarly, we have $[U_\omega,X_\omega]=(1+\omega(X_\Gamma))^{-1}U_\omega$.
These equations imply
 that $\Phi_\omega$ admits the $C^\omega$ Anosov splitting
 $TM_\Gamma=T\Phi_\omega \oplus \RR S_\omega \oplus \RR U_\omega$.

By Proposition \ref{prop:construction}
 there exists a $C^\omega$ Anosov flow $\ch{\Phi}_\omega$
 on a $C^\omega$ manifold $M_\omega$,
 a homeomorphism $H:M_\Gamma \ra M_\omega$,
 a $C^\omega$ vector field $Y^s$ on $M_\omega$,
 and a constant $\lambda>0$ such that
\begin{itemize}
\item $\ch{\Phi}^t_\omega \circ H=H \circ \Phi_\omega^t$
 and $D\ch{\Phi}_\omega^t(Y^s)=e^{-\lambda t} \cdot Y^s$ for any $t \in \RR$.
\item $J^u(H(p);\ch{\Phi}_\omega)=J^u(p;\Phi_\omega)$ for any
 $p \in \Per(\Phi_\omega)$.
\end{itemize}
Remark that
\begin{equation*}
 -J^s(H(p);\ch{\Phi}_\omega)
 =\lambda \cdot \tau(H(p);\ch{\Phi}_\omega)
 =\lambda \cdot \tau(p;\Phi_\omega)
 =\lambda  \cdot 
 \{ \tau(p;\Phi_\Gamma)+\Pair{\omega}{\gamma(p;\Phi_\Gamma)}\}
\end{equation*}
 for any $p \in \Per(\Phi_\omega)$.
Since $\Phi_\omega$ is a time-change of $\Phi_\Gamma$,
 we also have
 $\Per(\ch{\Phi}_\omega)=H(\Per(\Phi_\omega))=H(\Per(\Phi_\Gamma))$ and
\begin{equation*}
 J^u(H(p);\ch{\Phi}_\omega)=J^u(p;\Phi_\Gamma)=\tau(p;\Phi_\Gamma) 
\end{equation*}
 for any $p \in \Per(\Phi_\Gamma)$.

Let $\Phi_{\ch{\rho}_\omega}$ be the flow defined by
 $\Phi_{\ch{\rho}_\omega}^t=\ch{\Phi}_\omega^{t/\lambda}$
 and $\Psi_\omega$ the $C^\omega$ flow generated by $Y^s$.
Then, $\Phi_{\ch{\rho}_\omega}^t \circ \ch{\Psi}_\omega^x
 =\ch{\Psi}_\omega^{e^{-t}x} \circ \Phi_{\ch{\rho}_\omega}^t$
 for any $x,t \in \RR$.
Since $\Phi_{\ch{\rho}_\omega}$ is a constant time-change of $\ch{\Phi}_\omega$,
 we have $\tau(q;\Phi_{\ch{\rho}_\omega})
 =\lambda \cdot \tau(q;\ch{\Phi}_\omega)$,
 $J^u(q;\Phi_{\ch{\rho}_\omega}) = J^u(q;\ch{\Phi}_\omega)$,
 $J^s(q;\Phi_{\ch{\rho}_\omega}) = J^s(q;\ch{\Phi}_\omega)$
 for any $q \in \Per(\ch{\Phi}_\omega)$.
Hence,
\begin{align*}
J^u(H(p);\Phi_{\ch{\rho}_\omega}) 
 & = \tau(p;\Phi_\Gamma),\\
-J^s(H(p);\Phi_{\ch{\rho}_\omega}) 
 & = \tau(H(p);\Phi_{\ch{\rho}_\omega}) 
 = \lambda \cdot \left\{\tau(p;\Phi_\Gamma)
 +\Pair{\omega}{\gamma(p;\Phi_\Gamma)}\right\}
\end{align*}
 for any $p \in \Per(\Phi_\Gamma)$.
Now, a locally free $\GA$ action
 $\ch{\rho}_\omega$ on $M_\omega$ defined by
$\ch{\rho}_\omega(\ch{p},X^tS^x)
=\ch{\Psi}_\omega^x \circ \Phi_{\ch{\rho}_\omega}^{t}(\ch{p})$
 satisfies the required conditions.
\end{proof}

\begin{prop}
\label{prop:action}
If  $[\omega] \neq 0$,
 then the action $\ch{\rho}_\omega$ is not homogeneous.
\end{prop}
\begin{proof}
Let $H$ and $\lambda$ be the homeomorphism
 and the constant in Theorem \ref{thm:action}.
The equations in Theorem \ref{thm:action} imply
\begin{equation*}
J(H(p);\Phi_{\ch{\rho}_\omega}) 
 = (1-\lambda) \cdot \tau(p;\Phi_\Gamma)-\lambda \cdot
  \Pair{\omega}{\gamma(p;\Phi_\Gamma)}.
\end{equation*}

If $[\omega] \neq 0$,
 then there exists $p_1,p_2 \in \Per(\Phi_\Gamma)$ such that
 $\tau(p_1,\Phi_\Gamma)=\tau(p_2,\Phi_\Gamma)$ and
$\Pair{\omega}{\gamma(p_1,\Phi_\Gamma)}=
 -\Pair{\omega}{\gamma(p_2;\Phi_\Gamma)} \neq 0$.
Hence, at least one of
 $J(H(p_1);\Phi_{\ch{\rho}_\omega})$ or $J(H(p_2);\Phi_{\ch{\rho}_\omega})$
 is non-zero.
It implies that $\ch{\rho_\omega}$ admits no invariant volume.
Hence, $\ch{\rho}_\omega$ is not homogeneous.
\end{proof}

It is known that any homeomorphism between
 $C^\omega$ closed three-dimensional manifolds
 is homotopic to a $C^\omega$ diffeomorphism.
Therefore, we obtain
\begin{cor}
\label{cor:action}
If $H^1(M_\Gamma,\RR)$ is non-trivial,
 there exists a $C^\omega$ non-homogeneous locally free $\GA$-action
 on $M_\Gamma$.
\end{cor}

\section{Classification of actions}
\label{sec:classification}
In this section, we prove the main theorem.
After we review some known results on conjugacy between Anosov flows
 in Subsection \ref{sec:conjugacy},
 we introduce a natural map $\bar{a}_\Gamma:\AGamma \ra H^1(M_\Gamma,\RR)$
  in Subsection \ref{sec:bar-a}.
 It was originally introduced by Ghys \cite{Gh}
 as an obstruction to being a homogeneous action.
 We will see that the map $\bar{a}_\Gamma$
 classifies actions up to $C^\infty$ conjugacy homotopic to the identity.
In Subsection \ref{sec:Delta},
 we determine the image $\Delta_\Gamma$ of $\bar{a}_\Gamma$.
In fact, the actions constructed in Section \ref{sec:non-homogeneous}
 induce a family $(\rho_a)_{a \in \Delta_\Gamma}$ in $\AGamma$
 such that $\bar{a}_\Gamma(\rho_a)=a$.

%
%
\subsection{Conjugacy between Anosov flows}
\label{sec:conjugacy}
In this subsection, we review some known results
 on conjugacy between Anosov flows
 and give a criterion for $C^\infty$ conjugacy between
 locally free $\GA$-actions.

Let $\Phi_1$ and $\Phi_2$ be flows on manifolds $M_1$ and $M_2$,
 respectively.
We say that a homeomorphism $H:M_1 \ra M_2$
 is {\it a topological conjugacy}
 if $H \circ \Phi_1^t =\Phi_2^t \circ H$
 for any $t \in \RR$.
When $H$ is a $C^r$ diffeomorphism, it is called {\it a $C^r$ conjugacy}.
A continuous function $\alpha:M \times \RR \ra \RR$
 is called a {\it cocycle} over $\Phi_1$
 if it satisfies
 $\alpha(p,0)=0$ and
 $\alpha(p,s+t)=\alpha(p,s)+\alpha(\Phi_1^s(p),t)$
 for any $p \in M_1$ and $s,t \in \RR$.
\begin{Livthm}
[\cite{Li}]
\label{thm:Livschitz}
Let $\Phi$ be a $C^2$
 topologically transitive Anosov flow on a closed manifold $M$
 and $\alpha$ be a H\"older continuous cocycle over $\Phi$.
If $\alpha(p,\tau(p;\Phi))=0$ for any $p \in \Per(\Phi)$,
 then there exists a H\"older continuous function $\beta$ such that
\begin{equation*}
\alpha(p,t)=\beta \circ \Phi^t(p) -\beta(p)
\end{equation*}
 for any $p \in M$ and $t \in \RR$.
It is unique up to the constant term.
\end{Livthm}

The following results are applications of the Livschitz theorem.
\begin{thm}
[{\it c.f.} {\cite[Theorem 19.2.9]{KH}}]
\label{thm:Holder conjugacy}
Let $\Phi_1$ and $\Phi_2$ be $C^\infty$ topologically transitive
 Anosov flows on closed manifolds $M_1$ and $M_2$, respectively. 
Suppose that a topological equivalence $H$
 between $\Phi_1$ and $\Phi_2$ satisfies
 $\tau(p;\Phi_1)=\tau(H(p);\Phi_2)$ for any $p \in \Per(\Phi_1)$.
Then, there exists a H\"older continuous function $\beta$ on $M_1$
 such that the map $H_1:M_1 \ra M_2$ defined by
 $H_1(p)=\Phi_2^{\beta(p)}(H(p))$ is
 a topological conjugacy between $\Phi_1$ and $\Phi_2$.
\end{thm}
\begin{proof}
As Theorem 19.1.5 of \cite{KH},
 we can replace a topological equivalence $H$ by a bi-H\"older one.
Then, there exists a H\"older continuous cocycle
 $\alpha:M_1 \times \RR \ra \RR$ over $\Phi_1$ such that
 $H(\Phi_1^t(p))=\Phi_2^{\alpha(p,t)}(H(p))$
 for any $p \in M$ and $t \in \RR$.
The assumption implies
 $\alpha(p,\tau(p;\Phi_1))= \tau(H(p);\Phi_2)=\tau(p;\Phi_1)$
 for any $p \in \Phi_1$.
By the Livschitz theorem,
 there exists a H\"older continuous function $\beta$ on $M_1$
 such that $\alpha(p,t)=t+\beta(p)-\beta\circ \Phi_1(p)$.
Let  $H_1:M_1 \ra M_2$ be a continuous map defined
 by $H_1(p)=\Phi_2^{\beta(p)}(H(p))$.
We can show that $H_1 \circ \Phi_1^t=\Phi_2^t \circ H_1$
 by a direct computation.
It implies that $H_1$ is locally injective on each orbit of $\Phi_1$.
Hence, $H_1$ is a covering map.
Since the mapping degree of $H_1$ is one, it is a homeomorphism.
\end{proof}

We denote the topological entropy of a flow $\Psi$ by
 $h_{top}(\Psi)$.
\begin{prop}
\label{prop:SRB} 
Let $\Phi$ be a $C^2$ Anosov flow.
Suppose that there exists $\lambda>0$ such that
$J^u(p;\Phi)=\lambda \cdot \tau(p;\Phi)$ for any $p \in \Per(\Phi)$.
Then, $h_{top}(\Phi)=\lambda$.
\end{prop}
\begin{proof}
It follows from several standard facts
 in ergodic theory of hyperbolic systems.
Fix a Riemannian metric on $M$.
It is known that the subbundle $E^{uu}$ is H\"older continuous
 (see \cite[Theorem 19.1.6]{KH}).
By the Livschitz theorem,
 there exists a H\"older continuous function $\beta$ such that
 $\log|\det D\Phi^t|_{E^{uu}(p)}|
 =\lambda t+\beta \circ \Phi^t(p)-\beta(p)$
 for any $p \in M$ and $t \in \RR$.
Hence, the Lyapunov exponent along $E^{uu}$ is $\lambda$ at any point.
The entropy of $\Phi$-invariant measure
 is not greater than $\lambda$ by Ruelle's inequality,
 and the entropy of the SRB measure is $\lambda$ by Pesin's formula.
Hence, the variational principle implies that $h_{top}(\Phi)=\lambda$.
\end{proof}


De la Llave and ~Moriy\'on characterized
 $C^\infty$ conjugacy between three-dimensional Anosov flows
 by the coincidence of $J^s$ and $J^u$.
\begin{thm}
[\cite{LM2}]
\label{thm:Llave 1}
Let $\Phi_1$ and $\Phi_2$ be $C^\infty$
 topologically transitive Anosov flows
 on three-dimensional closed manifolds $M_1$ and $M_2$, respectively.
Suppose that
 a topological conjugacy $H$ between $\Phi_1$ and $\Phi_2$ satisfies
 $J^s(p;\Phi_1)=J^s(H(p);\Phi_2)$
 and $J^u(p;\Phi_1)=J^u(H(p);\Phi_2)$
 for any $p \in \Per(\Phi_1)$.
Then, $H$ is a $C^\infty$ diffeomorphism.
\end{thm}

As an application of the above results,
 we give a criterion of $C^\infty$ conjugacy between actions.
For a locally free $\GA$-action $\rho$
 on a closed three-dimensional manifold,
 we define a flow $\Phi_\rho$ by $\Phi_\rho^t(p)=\rho(p,X^t)$.
\begin{prop}
[{\cite[p.518]{Gh}}]
The flow $\Phi_\rho$ is topologically transitive.
\end{prop}
\begin{thm}
[\cite{As1}]	  
The flow $\Phi_\rho$ is Anosov.
\end{thm}
We call $\Phi_\rho$ {\it the Anosov flow associated with} $\rho$.

\begin{prop}
\label{prop:rho conj}
Let $\rho_1$ and $\rho_2$ be locally free $\GA$ actions
 on closed three-dimensional manifolds $M_1$ and $M_2$.
Let $\Phi_{\rho_i}$ be the Anosov flow associated with $\rho_i$.
Suppose that there exists a topological equivalence
 $H:M_1 \ra M_2$ between $\Phi_{\rho_1}$ and $\Phi_{\rho_2}$,
 and a constant $\lambda>0$ such that
$J^u(H(p);\Phi_{\rho_2})=J^u(p;\Phi_{\rho_1})$
 and $\tau(H(p);\Phi_{\rho_2})=\lambda \cdot \tau(p;\Phi_{\rho_1})$
 for any $p \in \Per(\Phi_{\rho_1})$.
Then, $\rho_1$ is $C^\infty$-conjugate to $\rho_2$ 
 by a diffeomorphism homotopic to $H$.
\end{prop}
\begin{proof}
First, we show that $\Phi_{\rho_1}$ and $\Phi_{\rho_2}$
 are topologically conjugate.
Put $\Psi^t=\Phi_{\rho_2}^{\lambda t}$.
Then, $H$ is a topological equivalence between $\Phi_{\rho_1}$ and $\Psi$.
Moreover,
 $J^u(H(p);\Psi)=J^u(p;\Phi_{\rho_1})$
 and $\tau(H(p);\Psi)=\tau(p;\Phi_{\rho_1})$
 for any $p \in \Per(\Phi_{\rho_1})$.
By Theorem \ref{thm:Holder conjugacy},
 there exists a H\"older continuous function
 $\beta$ on $M_1$ such that
 the map $H_1:M_1 \ra M_2$ defined by $H_1(p)=\Psi^{\beta(p)}(H(p))$
 is a topological conjugacy between $\Phi_{\rho_1}$ and $\Psi$.
In particular,
 $h_{top}(\Phi_{\rho_1})=h_{top}(\Psi)=\lambda \cdot h_{top}(\Phi_{\rho_2})$.
On the other hand,
 $h_{top}(\Phi_{\rho_1})=h_{top}(\Phi_{\rho_2})=1$ by Proposition \ref{prop:SRB}.
Therefore, $\lambda=1$, and hence, $\Psi=\Phi_{\rho_2}$.
It implies that
 $H_1$ is a topological conjugacy between $\Phi_{\rho_1}$ and $\Phi_{\rho_2}$,
 and $J^u(H_1(p);\Phi_{\rho_2})=J^u(p;\Phi_{\rho_1})$
 for any $p \in \Per(\Phi_{\rho_1})$.
 
Since $J^s(H_1(p);\Phi_{\rho_2})=-\tau(H_1(p);\Phi_{\rho_2})$
 and $J^s(p;\Phi_{\rho_1})=-\tau(p;\Phi_{\rho_1})$,
 we have $J^s(H_1(p);\Phi_{\rho_2})=J^s(p;\Phi_{\rho_1})$
 for any $p \in \Per(\Phi_{\rho_1})$.
By Theorem \ref{thm:Llave 1},
 the conjugacy $H_1$ is a $C^\infty$ diffeomorphism.

Let $X_i$ be the vector field generating $\Phi_{\rho_i}$
 and $S_i$ be the vector field generating the flow
 $(\rho_i(\cdot,S^x))_{x \in \RR}$ for each $i=1,2$.
Since $[S_i,X_i]=-S_i$,
 the vector field $S_i$ is tangent to
 the strong stable foliation of $\Phi_{\rho_i}$.
Since the $C^\infty$ conjugacy $H_1$
 preserves the strong stable foliation,
 there exists a $C^\infty$ function $g$ on $M_2$ such that
 $g(p) \neq 0$ and $DH_1(S_1)(p)=g(p) \cdot S_2(p)$ for any $p \in M_2$.
Then,
\begin{equation*}
g \cdot S_2=DH_1(S_1)=DH_1([X_1,S_1])
 =[X_2,g\cdot S_2]=(g+X_2 g) \cdot S_2,
\end{equation*}
 and hence, $X_2 g=0$.
Since the flow $\Phi_{\rho_2}$ is topologically transitive,
 the function $g$ is constant with a non-zero value $c$.
Now, we have
\begin{equation*}
 \rho_2(H_1(p),X^t S^{c x})=H_1(\rho_1(p,X^t S^x)) 
\end{equation*}
 for any $t,x \in \RR$ and $p \in M_\Gamma$.
It implies that $\rho_1$ is $C^\infty$-conjugate to $\rho_2$.
\end{proof}

\subsection{The map $\bar{a}_\Gamma$}
\label{sec:bar-a}
In this section, we introduce a natural parametrization of
 $\GA$-actions up to $C^\infty$ conjugacy homotopic to
 the identity map.

In \cite{Gh}, Ghys showed that
 the Jacobian of $\rho$ is controlled by a closed one-form.
\begin{thm}
[{\cite[Chapter IV]{Gh}}]
\label{thm:Ghys 2}
Let $M$ be a closed oriented three-dimensional manifold
 with a fixed volume
 and $\rho$ a $C^\infty$ locally free $\GA$-action on $M$.
Then, there exists a $C^\infty$ closed one-form $\omega_\rho$ on $M$
 such that
\begin{equation}
\label{eqn:def omega}
(\rho^g)^*\omega_\rho-\omega_\rho
 =-d\left(\log \det D\rho^g\right).
\end{equation}
 for any $g \in \GA$,
 where $\rho^g$ is a diffeomorphism on $M$ given by $\rho^g(p)=\rho(p,g)$.
\end{thm}

Fix a cocompact lattice $\Gamma$ of $\wSL(2,\RR)$
 and put $M_\Gamma=\Gamma \bsl \wSL(2,\RR)$.
The manifold $M_\Gamma$ is orientable
 and admits the standard volume
 induced from the Haar measure on $\wSL(2,\RR)$.
Let $\rho_\Gamma$ be the standard $\GA$-action on $M_\Gamma$
 and $\cF_\Gamma$ the orbit foliation of $\rho_\Gamma$.
Recall that the flow $\Phi_\Gamma$ on $M_\Gamma$ 
 is defined by $\Phi_\Gamma^t(\Gamma g)=\Gamma(g X^t)$.
It is the Anosov flow associated with $\rho_\Gamma$.
The following proposition characterizes
 an Anosov flow whose weak stable foliation is $\cF_\Gamma$.
\begin{prop}
\label{prop:isotopic equivalence}
Let $\Phi$ be an Anosov flow on $M_\Gamma$
 such that its weak stable foliation coincides with $\cF_\Gamma$.
Then, there exists a topological equivalence
 between $\Phi_\Gamma$ and $\Phi$
 which is isotopic to the identity and
 preserves each leaf of $\cF_\Gamma$.
\end{prop}
\begin{proof}
Let $\Sigma$ be a closed surface with a hyperbolic metric
 and $S^1\Sigma$ be its unit tangent bundle.
In Sections 3 and 4 of \cite{Gh6}, Ghys proved that
 if the weak stable foliation of an Anosov flow on
 $S^1 \Sigma$ is transverse to the fibers of
 the fibration $S^1 \Sigma \ra \Sigma$
 then there exists a topological equivalence
 between the Anosov flow and the geodesic flow on $S^1\Sigma$ 
 which is homotopic to the identity\footnote{Ghys did not
 mentioned that the topological equivalence he constructed
 is homotopic to the identity.
 However, this fact follows from his proof easily.}.
His proof works well even for our case
 by taking a finite covering of $S^1\Sigma$.
Hence, there exists a topological equivalence $H$ between
 $\Phi_\Gamma$ and $\Phi$ which is homotopic to the identity.

Let $\tilde{H}$, $\tilde{\Phi}$, $\tilde{\Phi}_\Gamma$
 be the lifts of $H$, $\Phi$, and $\Phi_\Gamma$ to $\wSL(2,\RR)$.
We denote the lift of $\cF_\Gamma$ to $\wSL(2,\RR)$ by $\tilde{\cF}_\Gamma$.
Since $\tilde{H}$ maps the $\tilde{\Phi}_\Gamma$-orbit of $x$
 to the $\tilde{\Phi}$-orbit of $\tilde{H}(x)$ for any $x \in \wSL(2,\RR)$,
 the Hausdorff distance between these orbits is bounded.
It implies that $\tilde{H}(x)$ is contained in $\tilde{\cF}_\Gamma(x)$.
Therefore, $H$ preserves each leaf of $\cF_\Gamma$.
\end{proof}

Recall that $\AGamma$ is the space of $C^\infty$ locally free
 $\GA$-action on $M_\Gamma$ whose orbit foliation is $\cF_\Gamma$.
By $\RR_+$, we denote the set of positive real numbers.
\begin{prop}
\label{prop:bar-a}
Let $\rho$ be an action in $\AGamma$,
 $H$ a topological equivalence
 between $\Phi_\Gamma$ and $\Phi_\rho$ which preserves
 each leaf of $\cF_\Gamma$,
  and $\omega_\rho$ a closed one-form given in Theorem \ref{thm:Ghys 2}.
Then, there exists a unique pair
 $(a,\lambda) \in H^1(M_\Gamma,\RR) \times \RR_+$
 such that
\begin{equation}
\label{eqn:bar-a}
\tau(H(p);\Phi_\rho)=
  \lambda \cdot  \{\tau(p;\Phi_\Gamma) +\Pair{a}{\gamma(p;\Phi_\Gamma)}\}
\end{equation}
 for any $p \in \Per(\Phi_\Gamma)$.
\end{prop}
\begin{proof}
Let $X_\rho$ be the vector field generating $\Phi_\rho$.
We define a function $\xi:M_\Gamma \ra \RR$
 by $\xi(p)=(d/dt)\log \det (D\Phi_\rho^t)_p|_{t=0}$.
By the definition (\ref{eqn:def omega}) of $\omega_\rho$,
\begin{align*}
 \omega_\rho(X_\rho) \circ \Phi_\rho^t-\omega_\rho(X_\rho)
 & = -X_\rho [\log \det D\Phi_\rho^t]\\
 &  = - \xi \circ \Phi_\rho^t(p)+\xi(p)
\end{align*}
 for any $p \in M_\Gamma$.
In particular, $\xi +\omega_\rho(X_\rho)$ is a $\Phi_\rho$-invariant function.
Since $\Phi_\rho$ is topologically transitive,
 there exists $\delta \in \RR$ such that $\xi+\omega_\rho(X_\rho)=\delta$.
It implies that
\begin{equation*}
 J(q;\Phi_\rho)
 = \delta \cdot \tau(q;\Phi_\rho)
  -\Pair{[\omega_\rho]}{\gamma(q;\Phi_\rho)}
\end{equation*}
 for any $q \in \Per(\Phi_\rho)$.
The map $H$ is homotopic to the identity map
 and preserves each leaf of $\cF_\Gamma$.
Hence,
 $[\gamma(H(p);\Phi_\rho)]=[\gamma(p;\Phi_\Gamma)]$ and
 $J^u(H(p);\Phi_\rho)=J^u(p;\Phi_\Gamma)=\tau(p;\Phi_\Gamma)$
 for any $p \in \Per(\Phi_\Gamma)$.
Since $\Phi_\rho$ is the Anosov flow associated with $\rho$,
 we have $J^s(H(p);\Phi_\rho)=-\tau(H(p);\Phi_\rho)$.
These equations imply
\begin{equation*}
 (1+\delta) \cdot \tau(H(p);\Phi_\rho)=
 \tau(p;\Phi_\Gamma) +\Pair{[\omega_\rho]}{\gamma(p;\Phi_\Gamma)}
\end{equation*}
 for any $p \in \Per(\Phi_\Gamma)$.
Since $[\gamma(p_0;\Phi_\Gamma)]=0$ for some $p_0 \in \Per(\Phi_\Gamma)$,
we have $1+\delta>0$.
Therefore, $[\omega_\rho]$ satisfies Equation (\ref{eqn:bar-a})
 for $\lambda=(1+\delta)^{-1}$.

We show the uniqueness of the pair $(a,\lambda)$.
Suppose that $(a,\lambda)$ satisfies the required condition.
Then,
\begin{equation*}
 \langle a-[\omega_\rho], \gamma(p;\Phi_\Gamma) \rangle
 = (\lambda^{-1}-(1+\delta))\tau(H(p);\Phi_\omega) 
\end{equation*}
 for any $p \in \Per(\Phi_\Gamma)$.
Take $p_0,\cdots,p_k \in \Per(\Phi_\Gamma)$ such that
 $[\gamma(p_0;\Phi_\Gamma)]=0$ and
 $[\gamma(p_1;\Phi_\Gamma)], \cdots,[\gamma(p_k;\Phi_\Gamma)]$ span
 $H_1(M_\Gamma,\RR)$.
The evaluation of the above equation
 at $[\gamma(p_0;\Phi_\Gamma)]$ implies $\lambda=(1+\delta)^{-1}$.
The evaluation at 
  $[\gamma(p_1;\Phi_\Gamma)], \cdots,[\gamma(p_k;\Phi_\Gamma)]$
 implies $a=[\omega_\rho]$.
\end{proof}

By the above proposition,
 the action $\rho$ determines the class $[\omega_\rho]$ uniquely.
We define maps $\bar{a}_\Gamma:\AGamma \ra H^1(M_\Gamma,\RR)$
 and $\bar{\lambda}_\Gamma:\AGamma \ra \RR_+$
 so that $(\bar{a}_\Gamma(\rho),\bar{\lambda}_\Gamma(\rho))$
 is the pair in Proposition \ref{prop:bar-a}.
\begin{prop}
\label{prop:injective} 
For $\rho_1,\rho_2 \in \AGamma$,
 $\bar{a}_\Gamma(\rho_1)=\bar{a}_\Gamma(\rho_2)$
 if and only if $\rho_1$ is $C^\infty$-conjugate to $\rho_2$
 by a diffeomorphism homotopic to the identity map.
\end{prop}
\begin{proof}
It is easy to see that $\bar{a}_\Gamma(\rho_1)=\bar{a}_\Gamma(\rho_2)$
 if two actions are $C^\infty$-conjugate
  by a diffeomorphism homotopic to the identity map.

Suppose that $\bar{a}_\Gamma(\rho_1)=\bar{a}_\Gamma(\rho_2)$.
By Proposition \ref{prop:isotopic equivalence},
 there exists a topological equivalence $H_i:M_\Gamma \ra M_\Gamma$ 
 between $\Phi_\Gamma$ and $\Phi_{\rho_i}$
 which is homotopic to the identity
 and which preserves each leaf of $\cF_\Gamma$.
Then, 
 $J^u(H_1(p);\Phi_{\rho_1}) =J^u(H_2(p);\Phi_{\rho_2}) =J^u(p;\Phi_\Gamma)$
 for any $p \in \Per(\Phi_\Gamma)$.
By Proposition \ref{prop:bar-a},
 the assumption $\bar{a}_\Gamma(\rho_1)=\bar{a}_\Gamma(\rho_2)$ implies
\begin{equation*}
 \bar{\lambda}_\Gamma(\rho_1)^{-1} \cdot \tau(H_1(p);\Phi_{\rho_1})
 =\bar{\lambda}_\Gamma(\rho_2)^{-1} \cdot \tau(H_2(p);\Phi_{\rho_2}).
\end{equation*}
Applying Proposition \ref{prop:rho conj} to $\rho_1$ and $\rho_2$,
 we obtain a diffeomorphism $H$ homotopic to $H_2 \circ H_1^{-1}$
 which conjugates $\rho_1$ to $\rho_2$.
Since $H_1$ and $H_2$ are homotopic to the identity map,
 so $H$ is.
\end{proof}

\subsection{The image of  $\bar{a}_\Gamma$}
\label{sec:Delta}
We define a subset $\Delta_\Gamma$ of $H^1(M_\Gamma,\RR)$ by
\begin{gather}
\label{eqn:def Delta} 
\Delta_\Gamma =
\left\{ a \in H^1(M_\Gamma,\RR) \;\left|\;
\sup_{p \in \Per(\Phi_\Gamma)}
 \frac{|\Pair{a}{\gamma(p;\Phi_\Gamma)}|}{\tau(p;\Phi_\Gamma)} <1
\right.
\right\}.
\end{gather}
We will show that the image of $\bar{a}_\Gamma$ is $\Delta_\Gamma$.

First, we show that
 each $a \in \Delta_\Gamma$ admits a nice representative.
\begin{lemma}
\label{lemma:ineq-Livschitz}
A cohomology class $a \in H^1(M_\Gamma,\RR)$
 is contained in $\Delta_\Gamma$
 if and only if
 it admits a $C^\omega$ representative
 such that $1+\omega(X_\Gamma)>0$.
\end{lemma}
\begin{proof}
The ``if'' part of the lemma is trivial.
We show the ``only if'' part.
By ${\cal M}(\Phi_\Gamma)$
 and ${\cal M}_{per}(\Phi_\Gamma)$,
 we denote the set of $\Phi_\Gamma$-invariant probability measures
 and its subset consisting of measures supported on a periodic orbit,
 respectively.

Take a $C^\omega$ closed one-form $\omega_0$ on $M_\Gamma$
 which represents $a \in \Delta_\Gamma$.
Then, there exists $\epsilon>0$ such that
\begin{equation*}
 \int 1+\omega(X_\Gamma) d\mu=
 1+\frac{\Pair{a}{\gamma(p;\Phi)}}{\tau(p;\Phi_\Gamma)} \geq \epsilon
\end{equation*}
 for any $\mu \in {\cal M}_{per}(\Phi_\Gamma)$.
It is known that ${\cal M}(\Phi_\Gamma)$ coincides with
 the convex hull of ${\cal M}_{per}(\Phi_\Gamma)$,
 see {\it e.g.} Lemma 2.4 of \cite{Gh1988}.
Hence, the integral of the function $(1+\gamma(X_\Gamma))$
 with respect to any $\mu \in {\cal M}(\Phi_\Gamma)$
 is positive.
By Lemma 2.5 of \cite{Gh1988},
 there exists a $C^\infty$ function $g_0$ on $M_\Gamma$ such that
 $1+\omega(X_\Gamma)+X_\Gamma g_0>0$.
Since $M_\Gamma$ is compact,
 we may approximate $g_0$ by a $C^\omega$ function $g$
 such that $1+\omega(X_\Gamma)+X_\Gamma g>0$.
The $C^\omega$ one-form $\omega+dg$
 satisfies the required condition.
\end{proof}

\begin{cor}
$\Delta_\Gamma$ is a convex open subset of $H^1(M_\Gamma,\RR)$.
\end{cor}

Next, we construct a family $(\rho_a)_{a \in \Delta_\Gamma}$ in $\AGamma$
 such that $\bar{a}_\Gamma(\rho_a)=a$ for any $a \in \Delta_\Gamma$.
Essentially, it is a corollary of
 the construction of non-homogeneous actions
 in Section \ref{sec:non-homogeneous}.
\begin{prop}
\label{prop:Delta 1} 
For any $a \in \Delta_\Gamma$,
 there exists $\rho_a \in \AGamma$ such that
 $\bar{a}_\Gamma(\rho_a)=a$.
\end{prop}
\begin{proof}
Fix $a \in \Delta_\Gamma$
 and a $C^\omega$ representative $\omega$ of $a$
 such that $1+\omega(X_\Gamma)>0$.
By Theorem \ref{thm:action},
 there exists a $C^\omega$ manifold $M_\omega$,
 a $C^\omega$ locally free $\GA$-action $\ch{\rho}_\omega$ on $M_\omega$,
 a topological equivalence $H$
 between $\Phi_\Gamma$ and $\Phi_{\ch{\rho}_\omega}$,
 and a constant $\lambda>0$ such that
\begin{align*}
 J^u(H(p);\Phi_{\ch{\rho}_\omega}) & = \tau(p;\Phi_\Gamma) ,\\
 \tau(H(p);\Phi_{\ch{\rho}_\omega})
 & = \lambda \cdot \left\{
 \tau(p;\Phi_\Gamma) +\Pair{[\omega]}{\gamma(p;\Phi_\Gamma)}
\right\}
\end{align*}
 for any $p \in \Per(\Phi_\Gamma)$.

Let $\cF_\omega$ be the orbit foliation of $\ch{\rho}_\omega$.
Since any locally free $\GA$-action on a closed oriented manifold 
 is $C^\infty$ orbit equivalent to a homogeneous action,
 there exists a homogeneous action $\rho_h$
 whose orbit foliation is $\cF_\omega$.
By Proposition \ref{prop:isotopic equivalence},
 $\Phi_{\ch{\rho}_\omega}$ is topologically equivalent
 to $\Phi_{\rho_h}$ by a homeomorphism $H_1:M_\omega \ra M_\omega$
 which is homotopic to the identity map
 and preserves each leaf of $\cF_\omega$.
Then, 
\begin{equation*}
\tau(H_1 \circ H(p);\Phi_{\rho_h})
 = J^u(H_1 \circ H(p);\Phi_{\rho_h}) 
 = J^u(H(p);\Phi_{\ch{\rho}_\omega})=\tau(p;\Phi_\Gamma)
\end{equation*}
 for any $p \in \Per(\Phi_\Gamma)$.
Since both $\rho_\Gamma$ and $\rho_h$ are homogeneous,
 it implies that
 $J^u(H_1 \circ H(p);\Phi_{\rho_h})=J^u(p;\Phi_{\rho_\Gamma})$
 and $J^s(H_1 \circ H(p);\Phi_{\rho_h})=J^s(p;\Phi_{\rho_\Gamma})$.
By Proposition \ref{prop:rho conj},
 there exists a $C^\infty$ diffeomorphism $H_\omega:M_\Gamma \ra M_\omega$
 which conjugates $\rho_\Gamma$ and $\rho_h$.

We define a $\GA$-action $\rho_a$ on $M_\Gamma$
 by $\rho_a(p,g)=H_\omega^{-1}(\ch{\rho}_\omega(H_\omega(p),g))$.
Since $H_\omega$ sends leaves of $\cF_\Gamma$ to those of $\cF_\omega$,
 $\rho_a$ is an element of $\AGamma$.
It is easy to check
 that $H_\omega^{-1} \circ H$ is a topological equivalence
 between $\Phi_\Gamma$ and $\Phi_{\rho_a}$ which is
 homotopic to the identity map and preserves each leaf of $\cF_\Gamma$.
Moreover, we have
\begin{align*}
\tau(H_\omega^{-1} \circ H(p);\Phi_{\rho_a})
& = \tau(H(p);\Phi_{\ch{\rho}_\omega})
 =\lambda \cdot \left\{
 \tau(p;\Phi_\Gamma)+\Pair{[\omega]}{\gamma(p;\Phi_\Gamma)}
 \right\}
\end{align*}
 for any $p \in \Per(\Phi_\Gamma)$.
Therefore, $\bar{a}(\rho_a)=[\omega]=a$.
\end{proof}

Combined with Propositions \ref{prop:injective} and \ref{prop:Delta 1},
 the following proposition completes the proof of the main theorem.
\begin{prop}
\label{prop:Delta 2} 
The image of $\bar{a}_\Gamma$ is contained in $\Delta_\Gamma$.
\end{prop}
\begin{proof}
Take $\rho \in \AGamma$.
Put $a=\bar{a}_\Gamma(\rho)$ and $\lambda=\bar{\lambda}_\Gamma(\rho)$.
Let $H$ be a topological equivalence between
 $\Phi_\Gamma$ and $\Phi_{\rho}$ which preserves
 each leaf of $\cF_\Gamma$.
Then, there exists $\epsilon>0$ such that
 $\epsilon \lambda \cdot \tau(p;\Phi_\Gamma) \leq \tau(H(p);\Phi_\rho)$
 for any $p \in \Per(\Phi_\Gamma)$.
Since $\tau(H(p);\Phi_\rho)=
 \lambda \cdot \{\tau(p;\Phi_\Gamma)+\Pair{a}{\gamma(p;\Phi_\Gamma)}\}$,
 we have
\begin{equation*}
\Pair{a}{\gamma(p;\Phi_\Gamma)}
 \geq -(1-\epsilon)\tau(p;\Phi_\Gamma)
\end{equation*}
 for any $p \in \Per(\Phi_\Gamma)$.

For any given $p \in \Per(\Phi_\Gamma)$,
 there exists $p' \in \Per(\Phi_\Gamma)$ such that
 $\tau(p';\Phi_\Gamma)=\tau(p;\Phi_\Gamma)$
 and $[\gamma(p';\Phi_\Gamma)]=-[\gamma(p;\Phi_\Gamma)]$.
The above inequality for $p'$ implies
$\Pair{a}{\gamma(p;\Phi_\Gamma)} \leq (1-\epsilon)\tau(p;\Phi_\Gamma)$.
Therefore, $a=\bar{a}_\Gamma(\rho) \in \Delta_\Gamma$.
\end{proof}


\appendix

%
%
\section{Locally free actions on solvable manifolds}
\label{section:solvable}
In this section,
 we give a proof of the following unpublished result by Ghys.
\begin{thm}
Let $M$ be a closed three-dimensional manifold
 whose fundamental group is solvable.
Then, any $C^\infty$ locally free action of $\GA$ on $M$ is homogeneous.
\end{thm}
\begin{proof}
Here, we present a shorter proof
 using a rigidity result due to Matsumoto and Mitsumatsu \cite{MM}.
Fix a $C^\infty$ locally free action $\rho$ on $M$.
Let $\cF_\rho$ be the orbit foliation of $\rho$.
By Plante's theorem \cite{Pl},
 the Anosov flow associated with $\rho$ is topologically equivalent to
 the suspension flow of a hyperbolic toral automorphism $A$.
Let $M_A$ be the mapping torus of $A$.
Since the orbit foliation of $\rho$ admits no closed leaves,
 a classification result by Ghys and Sergiescu \cite{GS}
 implies that $\cF_\rho$ is $C^\infty$ diffeomorphic to
 the orbit foliation $\cF_A$
 of a homogeneous $\GA$-action $\rho_A$ on $M_A$.
In \cite{MM},
 Matsumoto and Mitsumatsu proved that any locally free $\GA$-action
 on $M_A$ whose orbit foliation is $\cF_A$ is $C^\infty$-conjugate
 to $\rho_A$.
Therefore, $\rho$ is $C^\infty$-conjugate to the homogeneous action $\rho_A$.
\end{proof}

\section{Regularity of the unstable foliation}
In this section,
 we characterize homogeneous actions
 by the regularity of the unstable foliation  of the associated Anosov flow.
\begin{thm}
\label{thm:regularity} 
Let $\Gamma$ be a cocompact lattice of $\wSL(2,\RR)$
 and $\rho$ be an action in $\AGamma$.
If the weak unstable foliation of the Anosov flow
 associated with $\rho$ is a $C^2$ foliation,
 then $\rho$ is $C^\infty$ conjugate to
 the standard $\GA$-action on $M_{\Gamma}=\Gamma \bsl \wSL(2,\RR)$.
\end{thm}

For a cocompact lattice $\Gamma$ of $\wSL(2,\RR)$,
 let $\rho_{\Gamma}$ be the standard $\GA$-action on $M_\Gamma$
 and $\Phi_{\Gamma}$ the associated Anosov flow,
 {\it i.e.,} $\rho(\Gamma g,h)=\Gamma(gh)$
 and $\Phi_{\gamma}^t(\Gamma g)=\Gamma(g X^t)$.
For an Anosov flow $\Phi$,
 we denote the stable foliation and the unstable foliation of $\Phi$ by
 $\cF^s_\Phi$ and $\cF^u_\Phi$.
To simplify notations,
 we put $\cF^s_\Gamma=\cF^s_{\Phi_\Gamma}$
 and $\cF^u_\Gamma=\cF^u_{\Phi_\Gamma}$.
Remark that $\cF^s_\Gamma$ is the orbit foliation of
 the standard action $\rho_\Gamma$.
For an Anosov flow $\Phi$ on $M_\Gamma$ and $\gamma \in \Gamma$, we put
\begin{equation*}
\Per_\gamma(\Phi)=\{p \in M_\Gamma \st
 p=\Gamma g,\, 
 \gamma g =\tilde{\Phi}^T(g) \text{ for some }g \in \wSL(2,\RR),T>0\},
\end{equation*}
 where $\tilde{\Phi}$ is the lift of $\Phi$ to $\wSL(2,\RR)$.

The hyperbolic plane $\HH^2$
 admits a natural isometric action of $\wSL(2,\RR)$.
We define the {\it translation length} $L(g)$ of $g \in \wSL(2,\RR)$ by
\begin{equation*}
L(g)= \inf_{z \in \HH^2}d_{\HH^2}(z, g z),
\end{equation*}
 where $d_{\HH^2}$ is the distance induced by the hyperbolic metric.

Fix a cocompact lattice $\Gamma$ of $\wSL(2,\RR)$
 and an action $\rho \in \AGamma$.
Let $\Phi$ be the Anosov flow associated with $\rho$.
Suppose that $\cF^u_\Phi$ is a $C^2$ foliation.
\begin{lemma}
\label{lemma:regularity}
There exists an injective homomorphism $\sigma:\Gamma \ra \wSL(2,\RR)$
 such that
\begin{equation*}
L(\sigma(\gamma))= \bar{\lambda}_\Gamma(\rho)\left\{
L(\gamma)+\Pair{\bar{a}_\Gamma(\rho)}{\gamma}\right\}.
\end{equation*}
 for any $\gamma \in \Gamma$ with $\Per_\gamma(\Phi_\Gamma) \neq \emptyset$.
\end{lemma}
\begin{proof}
By a classification of three-dimensional Anosov flows
 with the $C^2$ Anosov splitting due to Ghys \cite{Gh4,Gh3},
 there exists a cocompact lattice $\Gamma'$ of $\wSL(2,\RR)$
 and a $C^2$ diffeomorphism $H:M_\Gamma \ra M_{\Gamma'}$ such that
 $H$ sends leaves of $\cF^u_\Phi$ to those of $\cF^u_{\Gamma'}$.
We define a flow $\Psi$ on $M_{\Gamma'}$
 by $\Psi^t=H \circ \Phi^t \circ H^{-1}$.
Remark that $\cF^u_\Psi$ coincides with $\cF^u_{\Gamma'}$.
The lift of $H$ to $\wSL(2,\RR)$ induces an isomorphism
 $\sigma:\Gamma \ra \Gamma'$ such that
 $H(\Per_\gamma(\Phi))=\Per_{\sigma(\gamma)}(\Psi)$
 for any $\gamma \in \Gamma$.

By Proposition \ref{prop:isotopic equivalence},
 there exists a topological equivalence $H_\Gamma:M_\Gamma \ra M_\Gamma$ 
 between $\Phi_\Gamma$ and $\Phi$ which preserves each leaf of $\cF^s_\Gamma$.
Similarly, there exists a topological equivalence
 $H_{\Gamma'}: M_{\Gamma'} \ra M_{\Gamma'}$
 between $\Phi_{\Gamma'}$ and $\Psi$
 which preserves each leaf of $\cF^u_{\Gamma'}$.

It is well-known and easy to see that
 $\tau(p;\Phi_\Gamma)=J^u(p;\Phi_\Gamma)=L(\gamma)$
 for any $\gamma \in \Gamma$ and $p \in \Per_\gamma(\Phi_\Gamma)$.
Similarly, $J^s(p';\Phi_{\Gamma'})=-L(\gamma')$
 for any $\gamma' \in \Gamma'$ and $p' \in \Per_{\gamma'}(\Phi_{\Gamma'})$.
Since $\gamma(H_{\Gamma'}(p');\Psi)$ is
 freely homotopic to $\gamma(p';\Phi_{\Gamma'})$ in $\cF^u_{\Gamma'}(p')$,
 we have $J^s(H_{\Gamma'}(p');\Psi)=J^s(p';\Phi_{\Gamma'}) =-L(\gamma')$
 for any $\gamma' \in \Gamma'$ and $p' \in \Per_{\gamma'}(\Phi_{\Gamma'})$.
Moreover, since $H$ is a $C^2$ diffeomorphism
 and $\Phi$ is the Anosov flow associated with $\rho$,
 we have
\begin{equation*}
 \tau(H_\Gamma(p);\Phi)=
 -J^s(H_\Gamma(p);\Phi)=-J^s(H \circ H_\Gamma(p);\Psi)=L(\sigma(\gamma))
\end{equation*}
 for any $\gamma \in \Gamma$ and $p \in \Per_\gamma(\Phi_\Gamma)$.
Now, the required equation follows
 from the definition of
 $(\bar{a}_\Gamma,\bar{\lambda}_\Gamma)$ and the above equations.
\end{proof}

Now we prove Theorem \ref{thm:regularity}.
Take $\gamma \in \Gamma$ with $\Per_\gamma(\Phi_\Gamma) \neq \emptyset$.
There exists $g \in \wSL(2,\RR)$ and $T>0$
 such that $\gamma g =g X^T$.
Put $g'=g \begin{bmatrix} 0 & -1 \\ 1 & 0 \end{bmatrix}$.
Then, $\gamma^{-1}g'=g'X^T$.
Hence, $\Per_{\gamma^{-1}}(\Phi_\Gamma)$ is non-empty.
Since $L(g^{-1})=L(g)$ for any $g \in \wSL(2,\RR)$,
 Lemma \ref{lemma:regularity} implies that
 $\Pair{\bar{a}_\Gamma(\rho)}{\gamma}=0$.
It is known that the set
\begin{equation*}
 \{[\gamma] \in H_1(M_\Gamma,\RR) \st
 \gamma \in \Gamma, \Per_\gamma(\Phi) \neq \emptyset\}
\end{equation*}
 spans $H_1(M_\Gamma,\RR)$.
Therefore, $\bar{a}_\Gamma(\rho)=0$.
It implies that $\rho$ is a homogeneous action.

%
%

\end{document}